\documentclass[xypic,amscd,syntonly,amssymb,verbatim,12pt]{amsart}

\usepackage{graphicx}
\usepackage[all]{xy}
\usepackage{amssymb}
\usepackage{amsmath}
\usepackage[mathscr]{euscript}
\usepackage{epsfig}
\theoremstyle{plain}
\newtheorem{Thm}{Theorem}
\newtheorem{Prop}[Thm]{Proposition}
\newtheorem{Cor}[Thm]{Corollary}
\newtheorem{Lem}[Thm]{Lemma}

 \theoremstyle{definition}
\theoremstyle{remark}

\numberwithin{equation}{section}

\begin{document}
 %\title{Coupling class of actions of reductive groups}
 \title{Equivariant branes}

 \author{ ANDR\'{E}S   VI\~{N}A}
\address{Departamento de F\'{i}sica. Universidad de Oviedo.   Avda Calvo
 Sotelo.     33007 Oviedo. Spain. }
 \email{vina@uniovi.es}
%\thanks{This work has been partially supported by Ministerio de Ciencia y
%Tecnolog\'{\i}a, grant FPA2009-11061}

 \maketitle
\begin{abstract}
%Considering the $B$-branes
 %on a Calabi-Yau manifold $X$ acted by a group $G$
 %as objects in the derived category of coherent sheaves, we give a definition of $G$-equivariant
%branes, which generalizes the concept of equivariant sheaves.
 Given a Calabi-Yau manifold $X$ acted by a group $G$ and considering
 the $B$-branes on $X$ as objects in the derived
 category of coherent sheaves, we give a definition of equivariant
branes, which generalizes the concept of equivariant sheaves.
We also propose a definition of equivariant charge of an equivariant
brane. The spaces of strings joining
 the branes ${\mathcal F}$ and ${\mathcal G}$, are the
  groups $Ext^i({\mathcal F},\,{\mathcal G})$.
 We prove that  the spaces of strings between two
$G$-equivariant branes  support representations of $G$.
 Thus, these spaces can be decomposed in
  direct sum of invariant spaces for the $G$-action. We show some
  particular decompositions, when $X$ is a toric variety and when $X$
  is a flag manifold of a semisimple Lie group.

\end{abstract}

\section {Introduction} \label{S:intro}

%In a first approximation, a
 As it is known, a
$D$-brane of type $B$ in a Calabi-Yau manifold $X$ is an object of
the derived category of coherent sheaves on $X$
\cite{Aspin,Aspin-et,Aspin-Law,Douglas,Katz,Sharpe}. In this note we will
consider   such objects in manifolds acted by a Lie group $G$.

Given a $G$-manifold $X$, some objects related  with $X$ admit an ``equivariant" version. For example,
 the equivariant vector bundles on $X$  are vector bundles equipped with a structure
compatible with the $G$-action on the base. In the same way,
% given a the manifold $X$  endowed with a $G$-structure,
 it is natural to consider ``equivariant"
$B$-branes on $X$. In this article, we deal with  equivariant branes on a $G$-manifold,
 with the spaces of open strings connecting them and we will
relate these spaces with representations of the group $G$.

Henceforth, the  space $X$ will be  a K\"ahler $G$-manifold, and
we put ${\mathcal O}_X$ for the corresponding structure sheaf. By
${\mathbf D}({\mathcal O}_X)$ we denote the bounded derived
category of coherent ${\mathcal O}_X$-modules
%over $X$
\cite{Kas-Sch}. In the context mentioned above, given
 the $B$-branes  ${\mathcal F}$ and ${\mathcal G}$,
that are objects of ${\mathbf D}({\mathcal O}_X)$, an open string
between  ${\mathcal F}$ and ${\mathcal G}$  is an element of the
Ext group ${Ext}^i_{{\mathcal O}_X}({\mathcal F},\,{\mathcal G})$
\cite{Aspin, Wei}, where $i+\,\hbox{ghost number of }\,{\mathcal
G}-\,\hbox{ghost number of }\,{\mathcal F}$ can be considered as
the ghost number of the corresponding strings.

 We denote by $\mu:G\times X\to X$ an {\em analytic} action of a reductive Lie group $G$ on $X$. Essentially,
 a $G$-equivariant structure on the ${\mathcal O}_X$-module  ${\mathcal H}$ is given by a family
 $\{\alpha_{g,x}\}$ of isomorphisms between the stalks
 $$\alpha_{g,x}:{\mathcal H}_x  \to{\mathcal H}_{\mu(g,x)},\;\;\;\hbox{for all}\;\;g\in G,\;x\in X$$
 compatible with the multiplication in $G$.

For a precise definition, we introduce the map $b:(g,x)\in G\times
X\mapsto x\in X$. A $G$-equivariant ${\mathcal O}_X$-module is a
pair $({\mathcal H},\,\alpha)$, where
 ${\mathcal H}$ is an ${\mathcal O}_X$-module  and $\alpha$ is  an isomorphism
\begin{equation}\label{b*mathcal}
\alpha:b^*{\mathcal H}\to \mu^*{\mathcal H}
 \end{equation}
 of ${\mathcal
O}_{G\times X}$-modules, where $b^*$ and $\mu^*$ are the functors
inverse image defined by the respective maps \cite[page
136]{Ge-Ma}. Furthermore, $\alpha$ must satisfy the cocycle condition
(see \cite[page 2]{Be-Lu} and equation (\ref{cocycle}) below).

The same definition of $G$-equivariance is applicable to an object
${\mathcal A}$ of the derived category ${\bf D}({\mathcal O}_X)$,
now  $b^*$ and $\mu^*$ are functors from  ${\bf D}({\mathcal
O}_X)$ to the derived category of ${\mathcal O}_{G\times
X}$-modules (see Subsection \ref{SusectCocycle}).

We will put $j:S\hookrightarrow X$ for the inclusion of an open
subset $S$ of $X$, and $j_!$  will denote the corresponding
functor direct image with compact support \cite[page
103]{Kas-Sch}.

In Section 2, we will prove the following theorems.
\begin{Thm}\label{Thmrepresentasheaf}
 If $({\mathcal H},\alpha)$ is a $G$-equivariant coherent sheaf on $X$ and $S$ is a $G$-invariant open subset of $X$,
 then the isomorphism $\alpha$ determines a representation of $G$ on $St^i(j_!({\mathcal O}_S),\,{\mathcal H})$, all $i$.
 In particular, each space $St^i({\mathcal O}_X,\,{\mathcal H})$
 carries a representation of $G$ induced by $\alpha$.
 \end{Thm}

\begin{Thm}\label{ThmrepresentaObj}
 If $({\mathcal G},\beta)$ and $({\mathcal F},\gamma)$ are $G$-equivariant  objects of the category ${\bf
 D}({\mathcal O}_X)$, then the isomorphisms $\beta $ and $\gamma$ determine a
 representation of $G$ on $St^i({\mathcal F},\,{\mathcal G})$, in a natural way.
 \end{Thm}

In the particular case that ${\mathcal F}$ and  ${\mathcal G}$ are
the sheaves of sections of $G$-equivariant vector bundles,
 %${\mathcal V}$ and ${\mathcal W}$,
  the action
  %whose existence assets Theorem of
   $g\in G$ on a morphism
 $$\Phi\in St^0\big({\mathcal F},\,{\mathcal G}\big)=
 Hom_{{\mathcal O}_X}\big({\mathcal F},\,{\mathcal G}\big)$$
    is given by the natural representation
 $(g\cdot\Phi)(\,-\,)=g\Phi(g^{-1}(\,-\,))$.
 % (see Example at the end of Section \ref{Sect.Equivariant}).

%ref{ThmSti} to objects in the derived category.

The decomposition of the representations of compact groups in
direct sum of irreducible representations
 %characterized by its character
  permits to classify the elements of a given space
$St^i({\mathcal F},\,{\mathcal G})$ in  subspaces, which
can be labelled by the characters of the corresponding irreducible
representations. Thus, we have the following theorem.

\begin{Thm}\label{Thmcharacters}
Let $G$ be a compact  group $G$. If  ${\mathcal F}$ and ${\mathcal
G}$ are $G$-equivariant objects of ${\bf D}(X)$
%  locally sheaves of ${\mathcal O}_X$-modules,
then for each $i$,
\begin{equation}\label{Sti}
St^i({\mathcal F},\,{\mathcal G})\simeq \bigoplus_{A}n_{A}
A,
 \end{equation}
 where
the sum runs over a  complete set of pairwise nonisomorphic representations  of $G$, $n_{A}$ is a natural number
 and $n_{A}A$ is the direct sum of $n_{A}$   summands of the
irreducible representation $A$.
%sub-representation of $St^i({\mathcal F},\,{\mathcalG})$ with character $\chi$.
\end{Thm}

%The results stated above will be proved en Section 2.
In Section \ref{Sect.Proofs}, we give   an equivariant version of the charge
 of a $G$-equivariant brane, that coincides with the usual one,
 when the group $G$ is trivial (see Subsection \ref{SubsectEqch}). This equivariant charge can be
 evaluated using the localization formulas in equivariant
 cohomology \cite{B-G-V}.

Examples manifolds in which a group action is an important
ingredient of its structure are the coadjoint orbits of a Lie
group and the toric manifolds. In Section 3, we will show the form
which  Theorem \ref{Thmcharacters} adopts in some examples of
pairs of ${\mathcal O}_X$-modules, when $X$ is a toric manifold
and when $X$ is a flag manifold of a semisimple group.

\smallskip

{\it Notations.} Besides the already introduced notations,  we
also use the following:
% conventions and notations:

The category of complex vector spaces will be denoted by
$\mathfrak{Vect}$, and  we let $D(\mathfrak{Vect})$ for
the corresponding bounded derived category.

Given a locally compact space $Z$,
% we will consider sheaves of
%complex vector spaces on $Z$ and we let $D(\mathfrak{Vect})$ for
%the derived category of complex vector spaces on $Z$.
if ${\mathcal R}$ a sheaf of ${\mathbb C}$-algebras on $Z$, the bounded derived
category of sheaves
%of complex vector spaces
 on $Z$ which are ${\mathcal R}$-modules is denoted by $D({\mathcal R})$.
 As usual, $\Gamma(Z,\, . \,)$ will be the functor global sections and
we put
 $$R\Gamma(Z,\,.\,):D( {\mathcal R}) \to D(\mathfrak{Vect})$$
 for its derived \cite{Kas-Sch,Wei}.
 The composition of this functor with the cohomology functor $H^i$
 is denoted by $R^i\Gamma(Z,\,.\,)$
 $$R^i\Gamma(Z,\,.\,):D( {\mathcal R}) \to \mathfrak{Vect}.$$

If $f:Y\to Z$ is a continuous map,
$$Rf_*:D(f^{*}{\mathcal R})\to D({\mathcal R})$$
 will denote the derived
functor of the direct image functor $f_*$.

 In
general, if $Z$ is a ringed space  the structure sheaf will
be denoted by ${\mathcal O}_Z$.

\smallskip

{\it Acknowledgement.} I thank to Diego Rodr\'{\i}guez-G\'omez for
pointing  the reference \cite{Aspin} out to me.

%%%%%%%%%%%%%%%%%%%%%%%%%%%%%%%%%%%%%%%%%%%%%%%%%%%%%%%%%%%%%%%%%%%%%%%%%%%%%%%%%%%%%%%%%%%%%%%%%%%%%
%%%%%%%%%%%%%%%%%%%%%%%%%%%%%%%%%%%%%%%%%%%%%%%%%%%%%%%%%%%%%%%%%%%%%%%%%%%%%%%%%%%%%%%%%%%%%%%%%%%%%%%

\section{Proofs of the results}\label{Sect.Proofs}

\subsection {The cocycle condition.} \label{SusectCocycle} To
formulate the cocycle condition, above mentioned, we introduce the
following notations
$$m:  G\times G\to  G,\;\;\;m(g_1,\,g_2)=g_1g_2.$$
$$p:G\times G\times X\to G\times X,\;\;\; p(g_1,\,g_2\,,x)=(g_2\,,x).$$
Thus, one has the maps $p,$ $ m\times 1_X$ and $1_G\times \mu$
from $G\times G\times X$ to $G\times X$ and the corresponding functors
\begin{equation}\label{functors}
\xymatrix{{\mathbf D}({\mathcal O}_X)\ar@/^/[r]^{b^*}  \ar@/_/[r]_{\mu^*}  & D({\mathcal O}_{G\times X})   \ar[r]^{p^*} \ar@/_1pc/[r]_{(1_G\times \mu)^*} \ar@/^2pc/[r]^{(m\times 1_X)^*} & D({\mathcal O}_{G\times G\times X})\,,
}
\end{equation}
where an asterisk as superscript is used
for denoting  the inverse image functor between the corresponding
derived categories.
The equalities
$$b\circ(m\times 1_X)=b\circ p,\,\; b\circ(1_G\times\mu)=\mu\circ p,\,\; \mu\circ(1_G\times\mu)=\mu\circ(m\times 1_X)$$
give rise to equalities between the respective compositions of the functors in (\ref{functors}).

Given an object  ${\mathcal A}$  of the category ${\mathbf
D}({\mathcal O}_X)$, an isomorphism $\alpha:b^*{\mathcal
A}\to\mu^*{\mathcal A}$
 %the isomorphism $\alpha$ in (\ref{b*mathcal})
 satisfies the cocycle condition if
\begin{equation}\label{cocycle}
(m\times 1_X)^*(\alpha)=(1_G\times\mu)^*(\alpha)\circ p^*(\alpha).
\end{equation}
In this case, we say that the pair $({\mathcal A},\,\alpha)$ is an
$G$-equivariant object.

Both sides of equation (\ref{cocycle}) are isomorphisms between
two objects on the derived category of ${\mathcal O}_{G\times
G\times X}$-modules;
 more precisely, between the objects $d_0^*{\mathcal A}$ and $d^*{\mathcal A}$, where
  $d,\,d_0: G\times G\times X \to X$ are defined by
$d_0(g_1,\,g_2,\,x)=x$ and $d(g_1,\,g_2,\,x)=(g_1g_2)x.$ In other
words, the cocycle condition means  the commutativity of the
following triangle
 \begin{equation}\label{triangleZ}
  \xymatrix{
 {\mathcal Z}_1\ar[dr]_{(m\times 1_X)^*(\alpha)}\ar[rr]^{p^*(\alpha)} &&{\mathcal Z}_{2}\ar[dl]^{(1_G\times \mu)^*(\alpha)}\\
 & {\mathcal Z}_{3}, }
  \end{equation}
where
$$\begin{aligned}\notag
&{\mathcal Z}_1:=p^*b^*({\mathcal A})=(m\times 1_X)^*b^*({\mathcal A}),\;\;{\mathcal Z}_2:=p^*\mu^*({\mathcal A})=(1_G\times\mu)^*b^*({\mathcal A}) \\ &{\mathcal Z}_3:=(m\times 1_X)^*\mu^*({\mathcal A})=(1_G\times\mu)^*\mu^*({\mathcal A}).
\end{aligned}$$
%$$
%\xymatrix{
% A \ar[dr]_{(m\times 1_X)^*(\alpha)}\ar[rr]^{p^*(\alpha)} && C
 %(\mu p)^*({\mathcal A})=(b(1_G\times\mu))^*({\mathcal A}) \ar[dl]^{(1_G\times \mu)^*(\alpha)}\\
% & B
%$$

%%%%%%%%%%%%%%%%%%%%%%%%%%%%%%%%%%%%%%%%%%%%%%%%%%%%%%%%%%%%%%%%%%%%%%%%%%%%%%%%%%%%%%%%%%%%%%%%%%%%%%%%%%
%%%%%%%%%%%%%%%%%%%%%%%%%%%%%%%%%%%%%%%%%%%%%%%%%%%%%%%%%%%%%%%%%%%%%%%%%%%%%%%%%%%%%%%%%%%%%%%%%%%%%%%

\smallskip
\subsection{Equivariant sheaves.}\label{SubsectEqSh}
Now we suppose that the $G$-equivariant object ${\mathcal A}$
is a coherent sheaf
${\mathcal H}$ on $X$. One can define the category ${\rm Coh}^G(X)$, whose objects are the $G$-equivariant coherent sheaves on $X$. If $({\mathcal H}',\,\alpha')$ and $({\mathcal H},\,\alpha)$ are objects in this category, a morphism in ${\rm Coh}^G(X)$ from $({\mathcal H}',\,\alpha')$ to $({\mathcal H},\,\alpha)$ is a sheaf morphism $f:{\mathcal H}'\to{\mathcal H}$ such that $\alpha b^*(f)=\mu^*(f)\alpha'$.
%, then, by the cocycle condition,  the following isomorphisms between
%the corresponding stalks of ${\mathcal H}$ form a commutative
%diagram
%$$
%\xymatrix{
% {\mathcal H}_x\ar[dr]_{(m\times 1)^*(\alpha)}\ar[rr]^{p^*(\alpha)} &&{\mathcal H}_{g_2x}\ar[dl]^{(1\times \mu)^*(\alpha)}\\
% & {\mathcal H}_{(g_1g_2)x}
%}
%$$

Given an open subset $U\subset X$ and $g\in G$, we put
$U_g:=\{g\}\times U\subset G\times X.$ If $({\mathcal H},\,\alpha)$ is an object of ${\rm Coh}^G(X)$,
 the restriction to $U_g$ of the morphism of sheaves
$\alpha$ is denoted  $\alpha|_{U_g}$
\begin{equation}\label{alpha(U}
\alpha|_{U_g}:b^*{\mathcal H}(U_g)={\mathcal H}(U)\longrightarrow
\mu^*{\mathcal H}(U_g)={\mathcal H}(gU).
\end{equation}

 Hence, one has the following proposition.
 \begin{Prop}\label{Propiso1}
  $\alpha|_{U_g}$ determines an isomorphism of complex vector spaces
 \begin{equation}\label{iso1}
 %\Gamma(U,\,{\mathcal H}) \stackrel{\sim}{\longrightarrow} \Gamma(gU,\,{\mathcal H}).
 {\mathcal H}(U) \stackrel{\sim}{\longrightarrow} {\mathcal H}(gU)
 \end{equation}
The image of $\sigma_U\in{\mathcal H}(U)$ will be denoted
$g\cdot\sigma_U$. Thus, the element $g\in G$ determines an
isomorphism
\begin{equation}\label{mathcalHx}
 \alpha_{g,x}:{\mathcal H}_x\to {\mathcal H}_{gx},
\end{equation}
 for any $x\in X$.
 \end{Prop}

\begin{Lem}\label{Lema}
Let $g,h$ be elements of $G$ and $U$ an open set of $X$, then
$$\alpha|_{(gU)_h}\circ\alpha|_{U_g}=\alpha|_{U_{hg}}.$$
\end{Lem}

{\it Proof.} We consider   the commutative triangle
(\ref{triangleZ}) when ${\mathcal A}$ is the sheaf ${\mathcal H}$
and we restrict this triangle to $\{h\}\times U_g\subset
G\times G\times X$. The restriction of $(m\times 1_X)^*(\alpha)$
is  the morphism
$${\mathcal Z}_1(\{h\}\times U_g )={\mathcal H}(U)\longrightarrow {\mathcal Z}_3(\{h\}\times U_g )={\mathcal
H}((hg)U)$$
induced by $\alpha$. Thus,  by (\ref{alpha(U}), the mentioned restriction is
 $\alpha|_{U_{hg}}$.

 The restriction of $p^*(\alpha)$ to  $\{h\}\times U_g$
 $${\mathcal Z}_1(\{h\}\times U_g )={\mathcal H}(U)\longrightarrow {\mathcal Z}_2(\{h\}\times U_g )={\mathcal
H}(gU)$$
is (\ref{alpha(U}).

Finally, we consider the restriction of $(1_G\times\mu)^*(\alpha)$. It is the morphism
 $${\mathcal Z}_2(\{h\}\times U_g )={\mathcal H}(gU)\longrightarrow {\mathcal Z}_2(\{h\}\times U_g )={\mathcal
H}(h(gU))$$
 induced by $\alpha$, and according to (\ref{alpha(U}) it is $\alpha|_{(gU)_h}$.
Then the lemma follows from   the commutativity of (\ref{triangleZ}).
 \qed

\smallskip

As a direct consequence of Lemma
\begin{equation}\label{alphahg}
\alpha_{h,gx}\circ\alpha_{g,x}= \alpha_{hg,x}
\end{equation}

 Since $gX=X$, $\alpha|_{X_g}$ is an automorphism of the vector space ${\mathcal H}(X)$,   from the Lemma, one deduces
 \begin{equation}\label{alphaxhg}
 \alpha|_{X_h}\circ\alpha|_{X_g}=\alpha|_{X_{hg}}.
 \end{equation}
 %That is,
 %$\rho_g:=\alpha|_{X_g}$ is an automorphism of the vector space $\Gamma(X,\,{\mathcal H})$.
 The following proposition is a consequence from (\ref{alphaxhg}).

 \begin{Prop}\label{CorRepr}
 The automorphisms
 $\{\rho_g:=\alpha|_{X_g}\}_g$
 %of $\Gamma(X,\,{\mathcal H})$
 form a representation of $G$ in the vector space ${\mathcal H}(X)$.
 %$\Gamma(X,\,{\mathcal H})$.
\end{Prop}

  \begin{Lem}\label{equivResolution}
  Given $({\mathcal H},\,\alpha)$ an object of ${\rm Coh}^G(X)$, there is an resolution of  $({\mathcal H},\,\alpha)$ in ${\rm Coh}^G(X)$
  consisting of injective ${\mathcal O}_X$-modules.
  \end{Lem}
{\it Proof.} Given a point $x\in X$, ${\mathcal H}_x$ is a
${\mathbb C}$-vector space, so is an injective ${\mathbb
Z}$-module (i.e. a divisible abelian group). We set
 ${I}_x:={\rm Hom}_{\mathbb Z}({\mathcal O}_x,\, {\mathcal H}_x)$, where ${\mathcal O}_x$ is
 the stalk of the sheaf
 ${\mathcal O}_X$ at the point $x$. Then the inclusion
 ${\mathcal H}_x\hookrightarrow { I}_x$ is an embedding  of the ${\mathcal O}_x$-module
 ${\mathcal H}_x$ in an injective
  ${\mathcal O}_x$-module (see  \cite[III.7]{M-L}, \cite[page 123]{R}).

 By means of the $I_x$ one can construct an injective ${\mathcal O}_X$-module ${\mathcal J}$
 in which ${\mathcal H}$ can be embedded.
 $${\mathcal J}=\prod_{x\in X}(j_x)_*(I_x),$$
 where $j_x:\{x\}\hookrightarrow X$ (see \cite[page 207]{Hart}).

 As it is well-known, ${\mathcal H}\overset{i}{\hookrightarrow}{\mathcal J}$ is the $0$th-term of an
 injective resolution of the ${\mathcal O}_X$-module ${\mathcal H}$.
 As ${\mathcal H}$ is $G$-equivariant, by (\ref{mathcalHx}), for each $g\in G$, there is an
 isomorphism of ${\mathbb Z}$-modules ${\mathcal H}_x\to{\mathcal H}_{gx}$. Since the $G$-action on $X$ is analytic, one has an isomorphism ${\mathcal O}_{gx}\to {\mathcal O}_x$. So, we have
 % turn determine
  an isomorphism $I_x\to I_{gx}$. Thus, for any open subset $U\subset X$,
  there is an isomorphism of ${\mathbb C}$-vector spaces ${\mathcal J}(U)\to{\mathcal J}(gU)$, making commutative the following diagram, where
   the horizontal sequences are exact
   \begin{equation}\label{DiagExact}
    \xymatrix{ 0\ar[r]  & {\mathcal H}(U) \ar[d]_{\alpha|_{U_g}}\ar[r]^{i_U} & {\mathcal J}(U)\ar@{-->}[d] \\
   0\ar[r] &  {\mathcal H}(gU)\ar[r]^{i_{gU}} & {\mathcal J}(gU)
  }
 \end{equation}
 The isomorphisms ${\mathcal J}(U)\to{\mathcal J}(gU)$
  give rise to an isomorphism of ${\mathcal O}_X$-modules,  $\alpha^0:b^*{\mathcal J}\to \mu^*{\mathcal J}$.

 From (\ref{alphahg}), it follows
% The cocycle condition  satisfied by $\alpha$ implies that
  the equality of
 ${\mathcal J}(U)\to{\mathcal J}(gU)\to {\mathcal J}(h(gU))$ and ${\mathcal J}(U)\to{\mathcal J}((hg)U).$
  That is, $\alpha^0$
 satisfies the cocycle condition. Hence, ${\mathcal J}$ is an object of ${\rm Coh}^G(X)$. By the commutativity of
  (\ref{DiagExact}), it follows that
  $i:{\mathcal H}\to{\mathcal J}$ a morphism in that category.

  Diagram (\ref{DiagExact}) can be continued with the cokernels of $i_U$ and $i_{gU}$.
  We denote
  %quotient maps   from ${\mathcal J}(V)\to {\mathcal C}(V)$, where
   ${\mathcal C}(V)={\mathcal J}(V)/{\mathcal H}(V)$ and put ${\mathcal C}^+$ for denoting the sheaf
   associated to the presheaf ${\mathcal C}$.
   There  are isomorphisms of ${\mathbb C}$-vector spaces induced  canonically  between
   the of vector spaces of the following {\em commutative} diagram.
  \begin{equation}\label{DiagExact1}
    \xymatrix{ 0\ar[r]  & {\mathcal H}(U) \ar[d]_{\alpha|_{U_g}}\ar[r]^{i_U} & {\mathcal J}(U)\ar[d] \ar[r] &{\mathcal C}(U)\ar[d]\ar[r]&{\mathcal C}^+(U)\ar[d] \\
   0\ar[r] &  {\mathcal H}(gU)\ar[r]^{i_{gU}} & {\mathcal J}(gU) \ar[r] & {\mathcal C}(gU) \ar[r] &  {\mathcal C}^+(gU)
  }
 \end{equation}

The term ${\mathcal J}^1$ of an injective resolution of ${\mathcal
H}$ can be obtained from ${\mathcal C}^+$, by embedding ${\mathcal
C}^+$ in an injective object, as we have made with ${\mathcal H}$. Hence, there exists an isomorphism $\alpha^1$ of ${\mathcal O}_{G\times X}$-modules,  making commutative the following diagram, where we put ${\mathcal J}^0$ for the preceding ${\mathcal J}$.
\begin{equation}\label{DiagExact12}
    \xymatrix{ 0\ar[r]  & b^*{\mathcal H} \ar[d]_{\alpha}\ar[r]^{b^*i} & b^*{\mathcal J}^0 \ar[d]_{{\alpha}^0} \ar[r]^{b^*{\partial}^0} & b^*{\mathcal J}^1\ar[d]^{{\alpha}^1}  \\
   0\ar[r] & \mu^* {\mathcal H}\ar[r]^{\mu^*i} & \mu^*{\mathcal J}^0 \ar[r]^{\mu^*{\partial}^0} & \mu^*{\mathcal J}^1 }
 \end{equation}
 Continuing the process, we obtain a complex ${\mathcal J}^{\bullet}$ which satisfies the Lemma.
 \qed

\begin{Prop}\label{Proprephcoho}
 If $({\mathcal H},\,\alpha)$ is a $G$-equivariant ${\mathcal O}_X$-module, then for
 each $i$ the cohomology group $H^i(X;\,{\mathcal H})$ supports a
 representation of $G$ induced by the isomorphism $\alpha$.
  \end{Prop}
{\it Proof.} Let ${\mathcal J}^{\bullet}$ the injective resolution
of ${\mathcal H}$ constructed in Lemma \ref{equivResolution}. As
${\mathcal J}^i$ is $G$-equivariant, by Proposition \ref{CorRepr},
the space the ${\mathcal J}^{i}(X)$ carries the representation
$\rho^i$ of $G$
 defined by $\rho^i_g=\alpha^i|X_g$.
 Since the diagrams
\begin{equation}\label{DiagExact2}
    \xymatrix{  {\mathcal J^i}(X) \ar[d]_{\rho^i_g}\ar[r]^{{d}^i} & {\mathcal J}^{i+1}(X)\ar[d]^{\rho_g^{i+1}} \\
   {\mathcal J}^i(X)\ar[r]^{{d}^i} & {\mathcal J}^{i+1}(X)
  }
 \end{equation}
  are commutative,
  one has a representation of $G$ on each cohomology group $h^i({\mathcal J}^{\bullet }(X))$ of
 the complex ${\mathcal J}^{\bullet}(X)$.
 As ${\mathcal J}^{\bullet}$ is an injective resolution of
  the ${\mathcal O}_X$-module ${\mathcal H}$ and
  $H^i(X,\,{\mathcal H})$ is by definition the cohomology group
  $h^i(\Gamma(X,\,{\mathcal J}^{\bullet }))$,  the proposition follows.
   \qed

\smallskip

 The arguments given in the proof of Proposition \ref{Proprephcoho} are also valid when $X$
 is substituted by a $G$-invariant open subset
  $S\stackrel{j}{\hookrightarrow}X$, thus $H^i(S,\,{\mathcal H})$ also carries a representation of $G$.
 %Let $U$ be an $G$-invariant open subset of $X$, and $j:U\hookrightarrow X$ the inclusion.

 \smallskip
 \noindent
 {\bf Proof of Theorem \ref{Thmrepresentasheaf}.} As we said, let
 $j_!$ for the corresponding functor direct image with compact
 support. The theorem follows from the above observation together with the following identites
 $$H^i(S,\,{\mathcal H})=R^i\Gamma(S,\,{\mathcal H})=Ext^i_{{\mathcal O}_X}(j_!({\mathcal O}_S),\,{\mathcal H}).$$
 \qed

Let $({\mathcal F},\,\gamma)$, $({\mathcal G},\,\beta)$ be
$G$-equivariant ${\mathcal O}_X$-modules.
 By Proposition \ref{CorRepr},   ${\mathcal F}(X)$ and
 ${\mathcal G}(X)$ support representations of $G$.
  %If $\tau$ is an element of either ${\mathcal F}(X)$ or
  %${\mathcal G}(X)$, then we set $g\tau$ by the section transformed under the action of $g$.
  We put ${\mathcal K}:={\mathcal Hom}_{\mathcal O_X}({\mathcal F},\,{\mathcal
G})$ for the sheaf of homomorphisms from ${\mathcal F}$ to
${\mathcal G}$. Given an open subset $U\subset X$, for $\Phi_U\in
{\mathcal K}(U)$ and $g\in G$, we define $g\cdot\Phi_U\in{\mathcal
K}(gU)$ as follows:

 Given $\sigma\in
{\mathcal F}(gU)$, then $g^{-1}\cdot\sigma\in {\mathcal F}(U)$,
with the notation introduced in Proposition \ref{Propiso1}. Then
$\Phi_U(g^{-1}\cdot\sigma)\in {\mathcal G}(U)$ and we put
 \begin{equation}\label{gcdotP}
 (g\cdot\Phi_U)(\sigma):=g\cdot(\Phi(g^{-1}\cdot\sigma)).
 \end{equation}
 So, we have constructed an isomorphism
 %denote by $\eta|_{U_g}$
 $$  \eta|_{U_g} :{\mathcal K}(U)\longrightarrow
{\mathcal K}(gU),\;\; \;\;  \Phi_U\mapsto g\cdot \Phi_U.$$
Moreover,
$$   \eta|_{(gU)_h}\circ \eta|_{U_g}=\eta|_{U_{hg}}.$$
 Therefore, the isomorphisms $\{\eta|_{X_g}\}_g$ define a representation
of $G$ on the space ${\mathcal K}(X)$. That is,
\begin{Prop}\label{ProprepExt}  Let $({\mathcal F},\,\gamma)$, $({\mathcal
 G},\,\beta)$ be
 $G$-equivariant ${\mathcal O}_X$-modules, then ${\mathcal K}(X)$, with ${\mathcal K}:={\mathcal Hom}_{{\mathcal O}_X}({\mathcal
 F},\,{\mathcal G})$, supports a representation of $G$ induced by the isomorphisms $\gamma$ and $\beta$.
  \end{Prop}

%%%%%%%%%%%%%%%%%%%%%%%%%%%%%%%%%%%%%%%%%%%%%%%%%%%%%%%%%%%%%%%%%%%%%%%%%%%%%%%%%%%%%%%%%%%%%%%%%%%%%%%%%%%%%%%%%%
%%%%%%%%%%%%%%%%%%%%%%%%%%%%%%%%%%%%%%%%%%%%%%%%%%%%%%%%%%%%%%%%%%%%%%%%%%%%%%%%%%%%%%%%%%%%%%%%%%

\subsection{Equivariant complexes.}\label{SubsectEqCom}
 The  results of Subsection \ref{SubsectEqSh} can be generalized to  case
 %corresponds to the case
 when ${\mathcal F}$ and ${\mathcal G}$ are $G$-equivariant objects of the category ${\mathbf
 D}({\mathcal O}_X)$, the derived category of coherent sheaves on $X$. As we said, a $G$-equivariant object
  of ${\mathbf D}({\mathcal O}_X)$ is a pair $({\mathcal A},\,\alpha)$ consisting
of an object ${\mathcal A}$ of ${\mathbf
 D}(X)$ and an isomorphism in the derived category of ${\mathcal O}_{G\times
 X}$-modules, which satisfies (\ref{cocycle}).

Given $g\in G$,  and open subset $U$ of $X$,
%and $G$-equivariant
%element $({\mathcal A},\,\alpha)$.
 we denote by $L_g$ the diffeomorphism
$$L_g:(g,\,x)\in U_g\mapsto (g,\,gx)\in (gU)_g.$$
Thus, $b\circ L_g=\mu:U_g\to gU$. As $L_g$ is a diffeomorphism the
functors $L^*_g$ and $(L_g^{-1})_*$ are the same,
%one can identify identifications of the following ${\mathbb C}$-vector spaces of sections.
%using the notations introduced above,
one has the following relations among objects of the
derived category $D(\mathfrak{Vect})$, when  $({\mathcal A},\,\alpha)$ is an equivariant object.
%, which are  correlative to the ones written before Proposition \ref{Propiso1}.
% We denote by $L_g$ the
%diffeomorphism
%$$L_g:(g,\,x)\in U_g\mapsto (g,\,gx)\in (gU)_g.$$
%Thus, $b\circ L_g=\mu:U_g\to gU$. As $L_g$ is a diffeomorphism the
%functors $L^*_g$ and $(L_g^{-1})_*$ are the same, one can identify
%identifications of
%the following ${\mathbb C}$-vector spaces of sections.
\begin{align}
  R\Gamma(U,\,{\mathcal A})=& R\Gamma(U_g, b^*{\mathcal A})\simeq
R\Gamma(U_g,\,\mu^*{\mathcal A})=R\Gamma(U_g,\,L_g^*b^*{\mathcal
A})  \notag \\
=& R\Gamma(U_g,\,R(L_g^{-1})_*b^*{\mathcal
A})=R\Gamma((gU)_g,\,b^*{\mathcal A})=R\Gamma(gU,\,{\mathcal
A}).\notag
\end{align}
That is, we have an isomorphism $R\Gamma(U,\,{\mathcal A})
\stackrel{\sim}{\longrightarrow} R\Gamma(gU,\,{\mathcal A}).$ In
particular, when $U=X$, for any $i$ there is an isomorphism
 $$R\Gamma(X,\,{\mathcal A})\stackrel{\Hat r_g}{\longrightarrow}
R\Gamma(X,\,{\mathcal A}).$$
 By the cocycle condition $\Hat r_{hg}=\Hat r_h\circ \Hat r_g$.

On the other hand, as $\alpha$ is an isomorphism between two
complexes,
 the $\alpha^i$'s intertwine with the
boundary operators; so, the representation $\Hat r$ induces a
representation $r$ on each space
 $R^i\Gamma(X,\,{\mathcal A})$.
%$$R^i\Gamma(X,\,{\mathcal A})\stackrel{r_g}{\longrightarrow}
%R^i\Gamma(X,\,{\mathcal A}).$$
% By the cocycle condition $r_{hg}=r_h\circ r_g$.
 Since
 $$ R^i\Gamma(X,\,{\mathcal A})={Ext}^i({\mathcal
O}_X,\,{\mathcal A})=H^i(X,\,{\mathcal A}),$$
 $H^i(X,\,{\mathcal A})$ carries a representation of $G$ induced by $\alpha$. This result is a
generalization of Proposition \ref{CorRepr}.

\smallskip
{\bf $G$-equivariant Horseshoe lemma.} The well-known Horseshoe
lemma (see \cite[page 349]{R}, \cite[page 37]{W}) admits a
$G$-equivariant version, which can be proved following the steps
of the proof for the no-equivariant case.   This equivariant
version can be formulated as follows:

 Let $({\mathcal
H'},\,\alpha')$,   $({\mathcal H},\,\alpha)$ and $({\mathcal
H''},\,\alpha'')$  be  $G$-equivariant ${\mathcal O}_X$-modules and
$$0\to{\mathcal H}' \to  {\mathcal
H}\to{\mathcal H}''\to 0$$
 be a short exact sequence in the category ${\rm Coh}^G(X).$
 If ${\mathcal H'}\to{\mathcal J}^{'\bullet}$ and ${\mathcal H''}\to{\mathcal
 J}^{''\bullet}$ are  resolutions of ${\mathcal H}'$ and ${\mathcal H}''$ (resp.) in ${\rm Coh}^G(X)$ consisting of injective
 ${\mathcal O}_X$-modules. Then there are a resolution ${\mathcal J}^{\bullet}$ of $({\mathcal H},\,\alpha)$ in ${\rm Coh}^G(X)$, formed by injective ${\mathcal O}_X$-modules, and morphisms between the resolutions such that
 $$0\to{\mathcal J'}\to {\mathcal J}\to{\mathcal J''} \to 0$$
 is an exact sequence of complexes in ${\rm Coh}^G(X)$.

\smallskip

 \noindent
 {\bf Proof of Theorem \ref{ThmrepresentaObj}.} Since $({\mathcal G},\beta)=\big(\{{\mathcal
G}^{\bullet}\},\,\{\beta^{\bullet}\}\big)$ is a $G$-equivariant
object of ${\mathbf D}(X)$, the ${\mathcal O}_X$-module $
{\mathcal G}^j$ together with $\beta^j$ is a $G$-equivariant
${\mathcal O}_X$-module. We denote by $\partial^j:{\mathcal G}^j\to{\mathcal G}^{j+1}$ the corresponding boundary operator, as $b^*$ and $\mu^*$ are exact functors, ${\rm ker}\,\partial^j$,  ${\rm im}\,\partial^j$ and the cohomology $h^j({\mathcal G}^{\bullet})$
% ${\rm ker}\,\partial^j/{\rm im}\,\partial^{j-1}$
 are $G$-equivariant ${\mathcal O}_X$-modules.

According to Lemma \ref{equivResolution}, there are resolutions in
${\rm Coh}^G(X)$ consisting of injective ${\mathcal O}_X$-modules
for  $h^j({\mathcal G}^{\bullet})$ and ${\rm im}\,\partial^{j-1}$.
% ${\rm ker}\,\partial^j/{\rm im}\,\partial^{j-1}$.
By  the $G$-equivariant  Horseshoe lemma applied to the exact
sequence
 $$0\to {\rm im}\,\partial^{j-1}\to {\rm ker}\,\partial^j\to h^j({\mathcal G}^{\bullet})\to 0,$$
 there is a resolution for ${\rm ker}\,\partial^j$
satisfying the properties above stated. A new application of the
equivariant Horseshoe lemma to the exact sequence
$$0\to{\rm ker}\,{\partial}^j\to{\mathcal G}^j\to{\rm im}\,{\partial}^j\to 0,$$
permits
%These constructions are the $G$-equivariant version of the first steps in
the construction of a $G$-equivariant Cartan-Eilenberg resolution ${\mathcal J}^{\bullet,\bullet}$ of the complex ${\mathcal G}^{\bullet}$, in which each term is an injective $G$-equivariant ${\mathcal O}_X$-module (see \cite[Theorem 10.45]{R}).
 %, according to the definition given at the beginning of Subsection \ref{SubsectEqSh}.

   The isomorphism between the complexes $b^*{\mathcal F}$ and  $\mu^*{\mathcal F}$ implies that the representations on ${\mathcal F}^a(X)$
   and on  ${\mathcal F}^{a+1}(X)$ satisfies $\partial^a(g\cdot \sigma)=g\cdot(\partial^{a}\sigma)$, for all $\sigma\in{\mathcal F}^a(X)$.
The total complex  ${\mathcal I}={\rm Tot}({\mathcal J^{\bullet,\bullet}})$ is a complex
 in ${\rm Coh}^G(X)$  quasi-isomorphic to ${\mathcal G}$ and it is formed by injective ${\mathcal O}_X$-modules.
Thus, by Proposition \ref{CorRepr},   each vector space ${\mathcal
I}^i(X)$ carries a representation of $G$, which also intertwine with
the boundary operators of ${\mathcal I}^{\bullet}(X)$.

  Since $({\mathcal F},\,{\gamma})$
  %=\{{\mathcal F}^{\bullet}\}$
  is a
  $G$-equivariant object of ${\mathbf D}(X)$, from the argument before Proposition \ref{ProprepExt}, it follows that
  ${\rm Hom}_{ {\mathcal O}_X(X)}\big({\mathcal
  F}^a(X),\,{\mathcal I}^b(X)\big)$ supports a representation
  of $G$. So,
  \begin{equation}\label{Prod}
  {\mathcal C}^n:=\prod_{a}{\rm Hom}_{ {\mathcal O}_X(X)}\big({\mathcal
  F}^a(X),\,{\mathcal I}^{a+n}(X)\big).
   \end{equation}
  carries also a representation $\rho$ of $G$.
  %, induced by $\gamma$ and $\beta$.
  For ${\mathcal C}^{\bullet}$ one defines the following boundary operator \cite[page 17]{Iversen}
  \begin{equation}\label{deltanfa}
  \delta^n(f_a)=\big(\partial^{a+n}\circ f_a+(-1)^{n+1}f_{a+1}\circ\partial ^a  \big),
  \end{equation}
  where,
 % On the other hand, given
    $$f_a\in{\rm Hom}_{{\mathcal O}_X(X)} \big({\mathcal
  F}^a(X),\,{\mathcal I}^{a+n}(X)\big).$$

  On the other hand, by (\ref{gcdotP}), the action of $g\in G$ on $f_a$ is
  given by  $g\cdot (f_a)(\sigma)=g\cdot f_a(g^{-1}\cdot \sigma), $ for all $\sigma\in{\mathcal F}^a(X)$.
  Since
  $$\partial^{a+n}\circ\big( g\cdot f_a(g^{-1}\sigma)  \big)=g\cdot\big( \partial^{a+n}(f_a(g^{-1}\cdot\sigma)  \big),\;\;\;
  \partial^a(g^{-1}\cdot\sigma)=g^{-1}\cdot\partial ^a(\sigma),$$
  then $\delta^n((g\cdot( f_a))=g\cdot \delta^n(f_a)$.
    Hence,  $\rho$  induces a representation $r$ of $G$ on the cohomology of ${\mathcal C}^{\bullet}$.

  On the other hand, the functor
  %Furthermore, the funcor(\ref{Prod})
  %is the $n$th term of the complex which determines
  $$R{\rm Hom}(\,.\,,\,.\,):{\mathbf D}(X)\times {\mathbf D}(X)\to D(\mathfrak{Vect})$$
  assigns to the pair $({\mathcal F},\,{\mathcal G})$ the object
  represented by the complex ${\mathcal C}^{\bullet}$. As
  % $n$th term is (\ref{Prod}).
  ${Ext}^i({\mathcal F},\,{\mathcal G})=h^i(R{\rm Hom}({\mathcal F},\,{\mathcal
  G})),$
   %the theorem follows,
   the representation  $r$ is the one claimed in the statement of the theorem.
 \qed

\smallskip

%%%%%%%%%%%%%%%%%%%%%%%%%%%%%%%%%%%%%%%%%%%%%%%%%%%%%%%%%%%%%%%%%%%%%%%%%%%%%%%%%%%%%%%%%%%%%%%%%%%%%%
%%%%%%%%%%%%%%%%%%%%%%%%%%%%%%%%%%%%%%%%%%%%%%%%%%%%%%%%%%%%%%%%%%%%%%%%%%%%%%%%%%%%%%%%%%%%%%%%%%%%%%%%%

 \subsection{ Equivariant charges.}\label{SubsectEqch}  The charge of a brane  ${\mathcal F}$ \cite{Aspin, Harvey, M-M} is an element
  of the cohomology of $X$ defined from certain characteristic clases of $X$ and ${\mathcal F}$.
  For a $G$-equivariant brane, it is natural to define a $G$-equivariant charge through the respective $G$-equivariant
  characteristic classes. The resulting charge will be an element of the equivariant cohomology  $H_G(X)$.

  %To each $G$-equivariant brane, we will associate an element of the equivariant cohomology of
 %$X$ (
 Given a $G$-equivariant brane ${\mathcal
 F}\in{\bf D}(X)$, the complex ${\mathcal F}^{\bullet}$ es
 quasi-isomorphic to a complex ${\mathcal E}^{\bullet}$ consisting of
 locally free sheaves; that is, ${\mathcal E}^i$ is the sheaf of
 sections of a $G$-equivariant  vector bundle ${\mathcal V}^i$.

Each ${\mathcal V}^i$ has the corresponding $G$-equivariant Chern
character ${\rm ch}^G({\mathcal V}^i)$ (see \cite[page
212]{B-G-V}). We put ${\rm ch}^G({\mathcal E}^{\bullet})$ for
denoting the $\sum_i (-1)^i{\rm ch}^G({\mathcal V}^i).$ On the
other hand, one can consider the $G$-equivariant $\Hat A$-class of
$X$, which will be denoted by $\Hat A^G(X)$. The equivariant
charge of the equivariant brane ${\mathcal  F}$ the can be defined
by the formula
 \begin{equation}\label{EquiCharge}
Q^G({\mathcal F}):= {\rm ch}^G({\mathcal E}^{\bullet})\Hat A^G(X).
 \end{equation}
 Taking into account the relation between the $\Hat A$-roof class
 and the Todd class \cite[page 231]{L-M} the above
 definition coincides, when   $G$ is the trivial group,
 with the one given in \cite{Aspin}.

 In some particular cases, the preceding definition has a natural interpretation
 in terms of the index  of an elliptic operator.
 The exterior bundle $\Lambda^*T^*X$  of $X$ with the connection
 induced by the Levi-Civita connection  and the standard
 Clifford multiplication  is a Dirac bundle (see \cite[page 114]{L-M}). This Dirac bundle
 has associated the corresponding Dirac operator $D^X$. If $G$ acts as a group of isometries
  of $X$, then $D^X$ is a $G$-operator \cite[page 211]{L-M}; i.e. $D^X$ is $G$-equivariant.

 Let us assume that the complex ${\mathcal E}^{\bullet}$ consists of only one
 nonzero element,  ${\mathcal E}^0$. The compactness of $G$ allows us to average over the group
 for obtaining $G$-invariant metrics and $G$-invariant connections on ${\mathcal V}^0$. On the other hand,
 the tensor product of $(\Lambda^*T^*X) \otimes{\mathcal V}^0$ is a Dirac bundle (see \cite[page 122]{L-M}) and the corresponding Dirac operator $D$ is also $G$-equivariant,  by the $G$-invariance of the metric and the connection.
  Since $D$ is elliptic, if $X$ is compact, ${\rm ker} \, D$ and ${\rm coker} \, D$ are representations
  of $G$ of finite dimension. Denoting by $R(G)$ the character ring of $G$,  the $G$-index of $D$, ${\rm ind}_G(D)$, is an element of $R(G)$. For $g\in G$ the virtual character $\chi(D)(g)$ of $D$ at $g$ is defined by
  $$\chi(D)(g)={\rm trace}(g|_{{\rm ker}\,D})-{\rm trace}(g|_{{\rm coker}\,D}).$$
  The equivariant index theorem \cite[Chapter 8]{B-G-V},
  %\cite{B-V}
  asserts that in a neighborhood of
   $0\in{\mathfrak g}:={\rm Lie}(G)$
  \begin{equation}\label{chi(D)}
   \chi(D)\circ{\rm exp}=\int_XQ^G({\mathcal E}^0).
   \end{equation}
 The value of $\chi(D)({\rm exp}(\xi))$, for $\xi\in{\mathfrak g}$,  can be calculated by the localization formula in
 equivariant cohomology. The result is the Atiyah-Segal-Singer
 fixed point formula \cite{B-G-V,L-M}.

 % is the equivariant index ${\rm ind}_G(D)$, for
 %the Dirac operator on  ${\mathcal V}^0$ associated to a
% $G$-invariant connection on  ${\mathcal V}^0$ \cite[Chapter 8]{G-G-K}. That
% is, for all $\xi$ in a neighborhood of $0\in{\mathfrak g}:={\rm
% Lie}(G)$
% $${\rm ind}_G(D)(\xi)=\int_X{\rm ch}^G({\mathcal V}^{0})(\xi)\Hat A^G(X)(\xi).$$

%%%%%%%%%%%%%%%%%%%%%%%%%%%%%%%%%%%%%%%%%%%%%%%%%%%%%%%%%%%%%%%%%%%%%%%%%%%%%%%%%%%%%%%%%%%%%%%%%%%%%%%
%%%%%%%%%%%%%%%%%%%%%%%%%%%%%%%%%%%%%%%%%%%%%%%%%%%%%%%%%%%%%%%%%%%%%%%%%%%%%%%%%%%%%%%%%%%%%%%%%%%%%%%%%%%

\section{Particular cases.}

In this section, we show the form that the results of Section 2
adopt in some simple cases.

 A consequence of the Grothendieck spectral
sequence \cite[page 207]{Ge-Ma} is the known Local-to-Global Ext
spectral sequence, which allows to determine the $Ext$ groups from
the sheaves ${\mathcal Ext}$. Given the ${\mathcal O}_X$-modules
${\mathcal F}$ and ${\mathcal G}$, the first quadrant spectral
sequence
$$E_2^{p,q}=H^p(X,\,{\mathcal Ext}^q({\mathcal F},\,{\mathcal G}))$$
abuts to $Ext^n({\mathcal F},\,{\mathcal G})$.

 Let ${\mathcal F}$ be a locally free ${\mathcal O}_X$-module of finite rank, then
% As ${\mathcal F}$ is locally free,
  $0\to{\mathcal F} \stackrel{1}{\to } {\mathcal F} \to 0$
 % a projective resolution
is a resolution  of ${\mathcal F}$,
 which can be used to determine the sheaves  ${\mathcal Ext}^q({\mathcal F},\,{\mathcal G})$
  \cite[Proposition 6.5, page 234]{Hart}. So, the sheaves ${\mathcal Ext}$ are the cohomology of the trivial complex consisting of the sheaf ${\mathcal Hom}({\mathcal F},\,{\mathcal G})$ at the position $0$ and zeros in the other positions. Thus,
 % $${\mathcal Ext}^0({\mathcal O}({\mathcal V}),\,{\mathcal O}({\mathcal W}))={\mathcal Hom} ({\mathcal O}({\mathcal V}),\,{\mathcal O}({\mathcal W}))$$
$${\mathcal Ext}^0({\mathcal F},\,{\mathcal G})= {\mathcal Hom}({\mathcal F},\,{\mathcal G}).$$
and ${\mathcal Ext}^q({\mathcal F},\,{\mathcal G})=0$, for $q\ne 0$. Therefore, the spectral sequence
 degenerates at the second page and
$$Ext^p({\mathcal F},\,{\mathcal G})=H^p(X,\, {\mathcal Hom}({\mathcal F},\,{\mathcal G})).$$
 %where ${\mathcal F}^{\vee}$ is the dual sheaf of ${\mathcal F}$.
 Thus, we have the following proposition.

\begin{Prop}\label{locallyfree}
Let ${\mathcal F}$ be a locally free sheaf of finite rank on $X$, then for any coherent sheaf ${\mathcal G}$,
$$St^q({\mathcal F},\,{\mathcal G})=H^q(X,\,{\mathcal F}^{\vee}\otimes {\mathcal G}),$$
where ${\mathcal F}^{\vee}$ is the dual sheaf of ${\mathcal F}$.
\end{Prop}
 Given an invertible sheaf ${\mathcal G}$, i.e. a rank $1$ locally free ${\mathcal
 O}_X$-module, we put ${\mathcal G}^n$ for denoting the tensor
 product ${\mathcal G}^{\otimes n}$. The following corollary is a consequence from the proposition together
  with the fact that ${\mathcal G}$ is ample (see  \cite[Proposition 5.3, page
  229]{Hart}).
\begin{Cor}\label{CorAmple}
Let ${\mathcal F}$ be a locally free sheaf of finite rank and
${\mathcal G}$ an ample invertible sheaf, then there is an integer
$n_0$ such that
$$St^q({\mathcal F},\,{\mathcal G}^n)=0,$$
for all $q>0$ and all $n>n_0$
\end{Cor}

 % For example, if  ${\mathcal V}$ is a holomorphic vector bundle on $X$ with finite rank, ${\mathcal F}$ can be the sheaf of holomorphic %sections of ${\mathcal V}$.

Let ${\mathcal V}$ and  ${\mathcal W}$  be holomorphic vector bundles over $X$ with finite rank.  We put ${\mathcal F}:={\mathcal O}(\mathcal V)$ and
 ${\mathcal G}:={\mathcal O}(\mathcal W)$ for denoting the respective sheaves of holomorphic sections. Then  ${\mathcal F}$ is a locally free sheaf and by the Proposition \ref{locallyfree},
 $$St^p({\mathcal O}({\mathcal V}),\,{\mathcal O}({\mathcal W}))=H^p(X,\, {\mathcal O}({\mathcal V}^{\vee}\otimes {\mathcal W})),$$
 where ${\mathcal V}^{\vee}$ is the dual vector bundle of ${\mathcal V}$.

Let us assume that ${\mathcal V}$ and ${\mathcal W}$ are $G$-equivariant vector bundles on $X$, with $G$ a compact Lie group. We denote by
  $\chi_V$ and $\chi_W$  the characters   of the corresponding representations on the spaces of sections. Then $St^0({\mathcal O}({\mathcal V}),\,{\mathcal O}({\mathcal W}))$ supports the representation with character $\chi:=\bar\chi_{V}\chi_W$. If we write
 $\chi=\sum_k n_k\chi_k$, where the $\chi_k$ are characters of a complete family of irreducible representations of $G$, then
 the open string states with ghost number $0$, between two  branes wrapped on the whole $X$ with gauge
 bundles ${\mathcal V}$ and
 ${\mathcal W}$, can be expressed as the following sum direct  of  $G$-invariant subspaces
 $$St^0({\mathcal O}({\mathcal V}),\,{\mathcal O}({\mathcal W}))\simeq\bigoplus_kn_kB_{k},$$
 where $B_k$ is a subspace on which the representation  of $G$ has character $\chi_k$.
 The natural number $n_k$ before the subspace $B_k$ means the direct sum of $n_k$ summands equal to $B_k$.

\smallskip

%%%%%%%%%%%%%%%%%%%%%%%%%%%%%%%%%%%%%%%%%%%%%%%%%%%%%%%%%%%%%%%%%%%%%%%%%%%%%%%%%%%%%%%%%%%%%%%%%%%%%%%%%%%%%%
%%%%%%%%%%%%%%%%%%%%%%%%%%%%%%%%%%%%%%%%%%%%%%%%%%%%%%%%%%%%%%%%%%%%%%%%%%%%%%%%%%%%%%%%%%%%%%%%%%%%%%%%%%%%%%%%%%%%%%

 \subsection{ Flag manifolds.}\label{SubsectFlag}
 When $X$ is a flag manifold of semisimple group, the result stated in Theorem
 \ref{Thmcharacters} admits, for
 certain spaces of strings, a more precise formulation derived
 from the Borel-Bott-Weil theorem (see \cite{B-S,W1}, for a brief exposition \cite[pages 13-22]{K-V}).
%  for certain spaces of strings.
 %We will consider the above results when $X$ is a flag manifold of
%semisimple group.

 We remind some basic facts about flag manifolds;
for details see \cite{Borel}. Let ${\mathfrak g}$ be the Lie
algebra of a linear connected semisimple complex Lie group
$G_{\mathbb C}$. We assume that  a Cartan subalgebra ${\mathfrak
h}$ of ${\mathfrak g}$ has been fixed.
% A Cartan subalgebra ${\mathfrak h}$ of ${\mathfrak g}$ is a
%maximal abelian subalgebra formed by semisimple elements. Let us assume that
%a Cartan subalgebra ${\mathfrak h}$ has been fixed.
 If ${\mathfrak v}$,  an ${\rm ad}({\mathfrak h})$-invariant
subspace of ${\mathfrak g}$, then the set roots of ${\mathfrak h}$
in ${\mathfrak v}$ will be denoted by  $\Delta({\mathfrak v})$.
Given a system ${\Delta}^+\subset\Delta({\mathfrak g})$ of
positive roots, the parabolic subalgebras of ${\mathfrak g}$
 can be constructed as
follows: If $\Gamma$ a set of simple positive roots, we put
$\Delta(\Gamma)={\rm Span}_{\mathbb
Z}(\Gamma)\cap\Delta({\mathfrak g})$. The set $\Gamma$ determines
the parabolic subalgebra ${\mathfrak p}:={\mathfrak
l}\oplus{\mathfrak u}$, where
 \begin{equation}\label{LeviDes}
 {\mathfrak l}:={\mathfrak
h}\oplus\bigoplus_{\alpha\in\Delta({\Gamma})}{\mathfrak
g}_{\alpha},\;\;\hbox{and}\;\;{\mathfrak u}:=
\bigoplus_{\alpha\in\Delta^+\setminus\Delta({\Gamma})}{\mathfrak
g}_{\alpha}.
 \end{equation}
 When ${\Gamma}=\emptyset$, the corresponding parabolic algebra is the  Borel subalgebra determined by ${\Delta}^+$.

  The parabolic subgroup $P$ associated
 with the algebra ${\mathfrak p}$ can be expressed $P=L_{\mathbb C}U_{\mathbb C}$ (Levi decomposition), where
 $L_{\mathbb C}\cap U_{\mathbb C}=\{1\}$,
 $L_{\mathbb C}$ is a reductive group and $U_{\mathbb C}$ is nilpotent.

 Henceforth in this Subsection, we assume that a parabolic subgroup $P$ of $G_{\mathbb C}$ has been fixed. As $G_{\mathbb C}$ is
  connected, the normalizer $N_{G_{\mathbb C}}({\mathfrak
 p})$
 %and $N_{G_{\mathbb C}}(P)$
 coincide with $P$. Hence, the  flag variety
 $X=G_{\mathbb C}/P$ can be identified with the set of all the parabolic subalgebras which are $G_{\mathbb C}$-conjugated to ${\mathfrak
 p}$.
%A well-known fact is that
 $X$ is a compact simply connected K\"ahler
 manifold (see for example \cite{Borel, W1}).  In \cite{Grant}, Grantcharov showed several examples
  flag manifolds which are Calabi-Yau.

By $G\subset G_{\mathbb C}$ we denote a real form of $G_{\mathbb
C}$; i.e. a Lie subgroup of $G_{\mathbb C}$, such that ${\mathfrak
g}= {\rm Lie}(G)\otimes_{\mathbb R}{\mathbb C}$.
 As subgroup of $G_{\mathbb C}$, $G$ acts on $X$ and there is only finitely
many $G$-orbits on $X$. In the case $G$ is
 a compact real form of $G_{\mathbb C}$, the  $G$-action on $X$ is transitive. In \cite{W} there is
  a detailed exposition of the properties $G$-action on $X$, a shorter one  can be looked up in \cite{Z}.

 We assume that $G$ is a {\em compact} real form of $G_{\mathbb C}$. We put $L:=L_{\mathbb C}\cap G$,
 then transitive action of $G$ on $X$ permits identify $X$ and $G/L$.
 Using the above Levi decomposition of $P$, an irreducible
representation $r$ of $L$ on the finite dimensional complex vector
space $V$ can be extended to a holomorphic representation of $P$,
in which the action of the factor $U_{\mathbb C}$ is trivial.
With the $P$-action on $V$ one can define the following
$G$-homogeneous vector bundle over $X=G_{\mathbb C}/P$
\begin{equation}\label{mathcalV}
{\mathcal V}=G\times_P V.
\end{equation}

Let $\lambda\in{\mathfrak h}^*$ be the highest weight of the above
representation $r$ of $L$. Define
\begin{equation}\label{ilambda}
 i({\lambda}):=\#\{\alpha\in\Delta^+\,|\,\langle
\lambda+\rho,\,\alpha\rangle<0  \},
 \end{equation}
  $\rho$ being the half sum of the positive roots.

 As $St^i({\mathcal O}_X,\,{\mathcal O}({\mathcal V}))=H^i(X,\,{\mathcal O}({\mathcal V}))$, by a direct application of  Borel-Bott-Weil theorem,
 we deduce the following proposition about the spaces of string on the flag manifold $X=G_{\mathbb C}/P$.
 \begin{Prop}\label{B-B-W}
  \indent
\begin{enumerate}
 \item If $\lambda+\rho$ is singular, that is $\langle\lambda+\rho,\,\alpha\rangle=0$ for some root
 $\alpha\in{\Delta}({\mathfrak g})$,  then $St^i({\mathcal
O}_X,\,{\mathcal O}({\mathcal V}))=0$ for all $i$.
 \item If $\lambda+\rho$ is regular,  that is $\langle\lambda+\rho,\,\alpha\rangle\ne
 0$ for all $\alpha\in{\Delta}({\mathfrak g})$,
 then there exists $w$ in the Weyl group so
 that $w(\lambda+\rho)$ is dominant with respect to
 $\Delta^+\cap\Delta({\mathfrak l})$. In this case,
 $St^i({\mathcal
O}_X,\,{\mathcal O}({\mathcal V}))=0$
 for $i\ne i(\lambda)$ and  $St^i({\mathcal
O}_X,\,{\mathcal O}({\mathcal V}))$ is the irreducible
representation of $G$ of highest weight
 $w(\lambda+\rho)-\rho$.
 \end{enumerate}
\end{Prop}

As a direct consequence of a vanishing theorem proved in
 \cite{S-W}, we can state the following result, which yields a upper bound for the ghost number of the strings
  between two particular types of  branes in $X$;  in other words,
  a upper bound in the number of nonzero summands of (\ref{Sti}).
\begin{Prop}\label{PropSchmWolf}
 Let $S$ be an open $G$-orbit ($G$ not necessarily compact) in the flag manifold $X=G_{\mathbb C}/P$. We denote by $j:S\hookrightarrow X$ the inclusion
  and by $s$ de complex dimension
 of a maximal compact subvariety of $S$. If  ${\mathcal H}$  is a coherent sheaf on
 $X$, then $St^i(j_!({\mathcal O}_S),\,{\mathcal H})=0$, for
 $i>s.$ Where $j_!({\mathcal O}_S)$ is the direct image of
 ${\mathcal O}_S$ by the inclusion.
 \end{Prop}

%%%%%%%%%%%%%%%%%%%%%%%%%%%%%%%%%%%%%%%%%%%%%%%%%%%%%%%%%%%%%%%%%%%%%%%%%%%%%%%%%%%%%%%%%%%%%%%%%%%%%%
%%%%%%%%%%%%%%%%%%%%%%%%%%%%%%%%%%%%%%%%%%%%%%%%%%%%%%%%%%%%%%%%%%%%%%%%%%%%%%%%%%%%%%%%%%%%%%%%%%%%%%%%%

\subsection{ Toric varieties.}\label{SubsectToric}
 Known properties of the cohomology of toric varieties will allow us to
 express the decomposition (\ref{Sti}) in a more precise terms, when
 ${\mathcal F}$ and ${\mathcal G}$ are particular branes on a
 toric manifold. We will also apply the localization formula in
 cohomology equivariant to (\ref{chi(D)}), when $X$ is a toric
 manifold.

Let $\Sigma$ be a fan in $N={\mathbb Z}^r$, we will denote by $X$
the toric variety defined by $\Sigma$ \cite{C-K,C-L-S,Fulton,Oda}.
We put $M:=Hom_{\mathbb Z}(N,\,{\mathbb Z}),$ $N_{\mathbb
R}:=N\otimes_{\mathbb Z}{\mathbb R}$ and $T$ for the torus
$$T=N\otimes{\mathbb C}^{\times}=Hom_{\mathbb Z}(M,\,{\mathbb
C}^{\times}).$$
  Given $m\in M$ we denote by $\chi^m$ the homorphism
$$\chi^m:t\in T\mapsto t(m)\in {\mathbb C}^{\times}.$$
 That is, the $\chi^m$'s are the characters of the irreducible representations of $T$.

 We put $\Sigma(1)$ for denoting
the set of $1$-dimensional cones in $\Sigma$, and given $\rho\in
\Sigma(1)$, there is a unique primitive element $v_{\rho}\in N\cap
\rho$, such that the cone $\rho$ can be expressed as ${\mathbb
R}_{\geq 0}v_{\rho}$, and any cone $\sigma$ of $\Sigma$ can be
written
 \begin{equation}\label{cone}
 \sigma=\sum_{\rho\in\Sigma(1)\cap\sigma}{\mathbb
R}_{\geq 0}v_{\rho}.
\end{equation}

We assume that $X$ is nonsingular. Given a family $(a_{\rho})\in
{\mathbb Z}^{\Sigma(1)}$, we put $\psi(v_{\rho})=a_{\rho}$, for
all $\rho\in\Sigma(1)$. By (\ref{cone}), $\psi$ can be extended to
%As the cone $\sigma$
 %$$\sigma=\sum_{\rho\in\Sigma(1)\cap\sigma}{\mathbb
%R}_{\geq 0}v_{\rho},$$
 %we have
  a  function $\psi$ defined on the
support of $\Sigma$, which is linear on each cone of $\Sigma$. On
the other hand, the family $(a_{\rho})$ determines the following
divisor on X
 \begin{equation}\label{Divisor}
 A=-\sum_{\rho}a_{\rho}V(\rho),
  \end{equation}
where $V({\rho})$ is the closure of the orbit of $\rho$ under the
$T$-action. $A$ is a $T$-invariant divisor of $X$, which
determines a $T$-equivariant line bundle ${\mathcal L}$, in the
usual way.

By Proposition \ref{Proprephcoho},
$H^i(X,\,{\mathcal O}({\mathcal L}))$ supports a representation of
$T$. As $Ext^i({\mathcal O}_X,\,{\mathcal O}({\mathcal L}))=
H^i(X,\,{\mathcal O}({\mathcal L}))$, the decomposition of
$St^i({\mathcal O}_X,\, {\mathcal O}({\mathcal
L}))$
%=H(X,\,{\mathcal O}({\mathcal L}))$
 stated in Theorem
\ref{Thmcharacters} can expressed, for this particular case, in
terms of the local cohomology of $N_{\mathbb R}$. In fact,   from
Theorem 2.6 of \cite{Oda} (see also \cite[page 74]{Fulton}) we
deduce the following proposition.
\begin{Prop}\label{PropToric}
If $\Sigma$ is a smooth fan and ${\mathcal L}$ is the line
$T$-equivariant bundle associated to the divisor (\ref{Divisor}),
then
$$St^i({\mathcal O}_X,\,{\mathcal O}({\mathcal L}))=\sum_{m\in
M}H^i_{Z(m)}(N_{\mathbb R},\,{\mathbb C})\chi^m,$$
 where $Z(m)=\{v\in N_{\mathbb R},|\, \langle m,\,v\rangle\geq
 \psi(v)\}$.
\end{Prop}

Let $A$ and $A'$ be $T$-invariant divisors of $X$. We denote by
$\psi$ and $\psi'$ the corresponding linear support functions
associated to $A$ and $A'$, respectively. We put ${\mathcal L}$
and ${\mathcal L}'$ for the respective line bundles. One says that
$\psi$ is strictly convex if
$\lambda\psi(u)+(1-\lambda)\psi(v)<\psi(\lambda u+(1-\lambda v)$,
for all $\lambda\in [0,\,1]$ and $u,v\in N_{\mathbb R}$. A known
fact is that ${\mathcal O}({\mathcal L})$ is ample if $\psi$ is
strictly convex \cite[page 70]{Fulton}. The following proposition
shows the form adopted by Theorem \ref{Thmcharacters}, when
${\mathcal F}={\mathcal O}({\mathcal L})$ and ${\mathcal
G}={\mathcal O}({\mathcal L'})$.
 \begin{Prop}\label{PropConvex}
 With the above notations, if $\psi'-\psi$ is strictly convex, then
 $$St^q({\mathcal O}({\mathcal L}),\,{\mathcal O}({\mathcal
 L}'))=0,\;\;\; \hbox{for}\;\;q\ne 0,$$
 and
  $$St^q({\mathcal O}({\mathcal L}),\,{\mathcal O}({\mathcal
 L}'))=\bigoplus_{m\in P} {\mathbb C}\cdot\chi^m,$$
 $P$ being $\{m\in M\,|\,\langle m,\, n\rangle\geq
 \psi'(n)-\psi(n),\;\hbox{for all}\; n\in N_{\mathbb R}\}$.
\end{Prop}
 {\it Proof.} By Proposition \ref{locallyfree},
 $St^q({\mathcal O}({\mathcal L}),\,{\mathcal O}({\mathcal
 L}'))=H^q(X,\,{\mathcal O}({\mathcal L}^{\vee}\otimes{\mathcal
 L'})).$ As the support function associated to ${\mathcal O}({\mathcal L}^{\vee}\otimes{\mathcal
 L'})$ is $\psi'-\psi$, the proposition follows from Theorem 2.7 and Corollary
 2.9 in \cite{Oda}.
  \qed

\smallskip

In order to apply the fixed point formula to (\ref{chi(D)}) when $X$ is a toric manifold, we make two Remarks.

{\it Remark 1.}
Let us assume that the torus $T$ acts {\em trivially} on a connected manifold   $S$ and that ${\mathcal W}$ is a $T$-equivariant vector bundle over $S$ with rank $m$.
By the equivariant splitting principle,
% for calculations,
we can assume that ${\mathcal W}$ is a direct sum of $T$-equivariant line bundles
 $${\mathcal W}=\bigoplus_{j=1}^m{\mathcal L}_j.$$
 The action of $T$ on ${\mathcal W}$ is defined by $m$ weights $\varphi_j$ and the $T$-equivariant
 Chern class of ${\mathcal L}_j$ is given by (see \cite[page 317]{G-G-K})
  \begin{equation}\label{cT1}
  c^T_1({\mathcal L_j})=c_1({\mathcal L}_j)+\frac{1}{2\pi}\varphi_j.
   \end{equation}
 The $T$-equivariant Chern character of ${\mathcal W}$ is (see \cite[page
  234]{L-M})
 \begin{equation}\label{chT}
 {\rm ch}^T({\mathcal W})=\sum_{j=1}^m{\rm exp}(c_1^T({\mathcal L}_j)).
  \end{equation}

{\it Remark 2.}
 Let $p$ be a fixed point for the $T$-action on the toric manifold associated with the fan $\Sigma\subset{\mathbb R}^n$.
 We denote by $\nu_{i,p}\in 2\pi({\mathbb Z})^n$,$\;i=1,\dots,n$,  the weights of the isotropy representation
 of $T$ on $T_pX$. The fixed points of the $T$-action are in bijective correspondence
with the $n$-dimensional cones in $\Sigma$ \cite[\textsection
3.2]{C-L-S}. If the point $p$ is associated with the cone
$\sigma$, then
  \begin{equation}\label{omegaip}
  \omega_{i,p}=\frac{\nu_{i,p}}{2\pi}
  \end{equation}
  are the generators of $\sigma^{\vee}\cap M$, where $\sigma^{\vee}$ is the dual cone of $\sigma$.

\smallskip

 In the statement of the following proposition, ${\mathcal V}^0$ is
 the vector bundle introduced in Subsection \ref{SubsectEqch}.
\begin{Prop}\label{PropFixedPoint}
Let $X$ be the toric manifold associated to the fan $\Sigma$, denoting by $\{\varphi_{j,p}\}_{j=1,\dots,m}$ the weights of the representation of $T$ on the fibre of ${\mathcal V}^0$ at $p$, then  (\ref{chi(D)}) is equal to
%$$\chi(D)\circ{\rm exp}=\frac{1}{2^n}\sum_{p\in X^T}\Big(\sum_{j=1}^m{\em e}^{\varphi_j,p} \Big)\Big(\sum_{i=1}^n\frac{1}{\sinh %\\nu_{i,p}}  pi^n
\begin{equation}\label{2-n}
 2^{-n}\sum_{p\in X^T}\Big(\sum_{j=1}^m{\rm
e}^{\frac{1}{2\pi}\varphi_{j,p}}
\Big)\prod_{i=1}^n\Big(\sinh\big(\frac{
\omega_{i,p}}{2}\big)\Big)^{-1},
 \end{equation}
 where $X^T$ is the set of fixed points of $X$ for the $T$-action.
\end{Prop}
{\it Proof.} %As the fixed point set for $T$-action is a discrete
%set,
The localization theorem  in equivariant cohomology
\cite[\textsection 10.9]{G-S} allows us to calculate the value
(\ref{chi(D)}) as a sum of contributions of the connected
components of $X^T$. As $X^T$ is discret, the localization formula
adopts the following form
%For $\xi$ in a neighborhood of ${\rm Lie}(T)$,
\begin{equation}\label{chi(D)1}
%\chi(D)({\rm exp}(\xi))=\sum_{p\in X^T}\frac{Q^G({\mathcal E}^0)(p)}{\prod_i\alpha_{i,p}}=\sum_{p\in X^T}\frac{m}{\prod_i\alpha_{i,p}},
 \chi(D)\circ{\rm exp}=(2\pi)^n\sum_{p\in X^T}\frac{Q^G({\mathcal E}^0)(p)}{\prod_i\nu_{i,p}},
 \end{equation}
where the $\nu_{i,p}$ are the weights of the isotropy representation of $T$ at the fixed point $p$.

From (\ref{cT1}) and (\ref{chT}), it follows
$${\rm ch}^T({\mathcal V}^0){\big|_p}=\sum_{j=1}^m{\rm e}^{\frac{1}{2\pi}\varphi_{j,p}}.$$
Similarly (see \cite[page 231]{L-M}),
$$\Hat A^T(TX){\big|_p}=\prod_{i=1}^n\frac{\omega_{i,p}}{2}\Big(\sinh\big(\frac{\omega_{i,p}}{2}\big)\Big)^{-1}.$$
The proposition follows from (\ref{EquiCharge}) together with
(\ref{omegaip}). \qed

Note that
%, as it is usual in toric geometry,
 the contribution of the manifold  $X$ to (\ref{2-n}) is encoded in the $n$-dimensional
cones of the fan $\Sigma$.

\end{document}